\documentclass[3p]{elsarticle}
\usepackage{graphicx}
\usepackage{amssymb}
\usepackage{amsthm}
\usepackage{lineno}

\usepackage{epstopdf}
\DeclareGraphicsExtensions{.eps}

\biboptions{sort&compress}

\newtheorem{thm}{Theorem}

\newdefinition{rmk}{Remark}
\newproof{pf}{Proof}
\newproof{pot}{Proof of Theorem \ref{thm2}}
\begin{document}
\begin{frontmatter}
\title{Application of the Parallel Dichotomy Algorithm for solving Toeplitz tridiagonal systems
of linear equations  with one right-hand side.}

\author{Andrew V. Terekhov}
\ead{andrew.terekhov@mail.ru}
\address{Institute of
Computational Mathematics and Mathematical Geophysics,
630090,Novosibirsk,Russia}
\address{Budker Institute of Nuclear Physics, 630090, Novosibirsk,
Russia}
\address{Novosibirsk State University, 630090, Novosibirsk,
Russia}

\begin{abstract}

Based on a modification of the "Dichotomy Algorithm"\ (Terekhov,
2010), we propose a parallel procedure for solving tridiagonal
systems of equations with Toeplitz matrices. Taking the structure
of the Toeplitz matrices, we may substantially reduce the number
of the "preliminary calculations" of the Dichotomy Algorithm,
which makes it possible to effectively solve a series as well as a
single system of equations. On the example of solving of elliptic
equations by the Separation Variable Method, we show that the
computation accuracy is comparable with the sequential version of
the Thomas method, and the dependence of the speedup on the number
of processors is almost linear. The proposed modification is aimed
at parallel realization of a broad class of numerical methods
including the inversion of Toeplitz and quasi-Toeplitz tridiagonal
matrices.
\end{abstract}
\begin{keyword}
Parallel Dichotomy algorithm \sep Tridiagonal matrix algorithm
(TDMA) \sep Thomas algorithm  \sep Poisson equation \sep Toeplitz
matrixes \sep Fourier Method

\PACS 02.60.Dc \sep 02.60.Cb \sep 02.70.Bf \sep 02.70.Hm
\end{keyword}

\end{frontmatter}

\section{Introduction}
The realization of many numerical algorithms, such as the
Multigrid Line Relaxation
\cite{Elliptic:Multi1,Elliptic:Multi2,Elliptic:Multi3}, ADI and
Separation Variable Method
\cite{Elliptic:Multi2,Hockney:Particles,Samarski_Nikolaev,Elliptic:Spectral1},
Cyclic Reduction Method \cite{Hockney:Particles,Samarski_Nikolaev}
for solving elliptic equations and also the problems of
construction of splines \cite{Boor1978,Shikin1995} etc. requires
solving tridiagonal systems of equations with Toeplitz matrices
\cite{Toeplitz,MatrixInvert}

\begin{equation}
\label{Toelitz_matrix} T=\left(
\begin{array}{ccccc}
  t_0 & t_1 &  &   &  \Large 0\\
  t_{-1} & t_0 & t_1  &  &  \\
%   & t_{-1} & t_0 & t_1  &  \\
    & \ddots & \ddots & \ddots &  \\
    &  & t_{-1} & t_{0} & t_{1} \\
   0&   &  & t_{-1} & t_0 \\
\end{array}
\right)\equiv\mathrm{tridiag}\{...,t_{-1},t_0,t_1,...\}.
\end{equation}

Since, in solving modern problems of mathematical modeling, the
size and the number of such problems can reach several tens of
thousands, such computations must be performed on a supercomputer.

A lot of versions of parallel Thomas algorithms for general
tridiagonal matrices
\cite{Tridiagonal2,Tridiagonal4,Tridiagonal5,Tridiagonal:Stone,tridiag-exmp1,Ortega:Parallel,Konovalov,Wang},
 as well as for Toeplitz matrices
\cite{Tridiagonal:Toeplitz2,Tridiagonal:Toeplitz,Tridiagonal:Toeplitz5},
 have been worked out as of today. Algorithms have been proposed
that take into account the presence of diagonal
dominance\cite{Tridiagonal3}. In \cite{Povitski:ADI},on the
example of the realization of the ADI method, the approach was
considered in which combination of computations and interprocessor
exchanges makes it possible to increase the paralleling
efficiency. In \cite{Terekhov:Dichotomy}, a method called the
Dichotomy Algorithm was proposed for solving a series of problems
with a constant matrix and different right-hand sides

\begin{equation}
\label{main_eq3} \begin{array}{c}
 A {\bf X_{m}}= { \bf F_{m}} ,\quad m=1,...,M,\\\\
 A=\mathrm{tridiag}\{...,c_i,b_i,a_i,...\},\;i=1,...N,\;c_1=a_N=0,
\end{array}
\end{equation}
where $M$ is the number of problems in the series. The advantage
of the Dichotomy Algorithm consists in that it makes it possible
to attain a thousandfold speedup for problem  (\ref{main_eq3}) not
only in theory but also in practice. Aiming at a further
development of this approach, we, based on a modification of the
Dichotomy Algorithm, will propose a parallel algorithm for solving
not only a series but a single system of linear algebraic
equations (SLAE) with Toeplitz matrices.

A sufficient condition for the applicability of the Dichotomy
Algorithm is the diagonal dominance of the matrix of the SLAE. As
regards accuracy, the number of arithmetic operations, and the
number of communication interactions, the Dichotomy Algorithm is
almost equivalent to the Cyclic Reduction method
\cite{Snonkwiler:Parallel,Hockney:Cyclic} However, for comparable
data volumes, the real time of interprocessor interactions in the
Dichotomy Algorithm is much smaller. This is explained by the fact
that all interprocessor exchanges may be carried out via a
sequence of calls to the collective operation
"All-to-One-Reduce(+)" \cite{Janssen}.The account of such
properties as the associativity and commutativity of the operation
"+" taken over the distributed data makes it possible to reduce
the time of communication interactions by optimizing
them\cite{Collective:1,Tuning:1,Tuning:2,Tuning:3,Optimizing3}.While
the organization of exchanges via a multiple call to the
nonblocking\footnote{This requirement is necessary for the
realization of the Dichotomy Algorithm (see Subsection
\ref{Par:cost})} function "All-to-One-Reduce(+)"\  makes it
possible to reduce the time of synchronization of processor
elements (PE). Indeed, if, in one group of processors \footnote{A
communicator in the terminology of the MPI.}, there exist two free
PEs with prepared data then the execution of the collective
operation "(+)" over this group of processors can start to even if
the previous call on all processors is not finished. Thus, the
structure of the Dichotomy Algorithm yields ample opportunities
for the minimization of the data transfer time as well as the PE
synchronization time. In the cyclic reduction method, the fixed
order of elimination of the unknowns, on the one hand, restricts
optimization of communication interactions, and on the other,
requires synchronization of the computations on each reduction
order.

The presence of a preliminary step with $O(N)$ arithmetic
operations spent, where $N$ is the dimension of the SLAE, does not
make it possible to use the Dichotomy Algorithm effectively for
solving one problem. However, for Toeplitz matrices, an economical
preliminary procedure can be constructed with the number of the
arithmetic operations of  order $O(N/p+\log_2 p)$, where $p$ where
is the number of the processors. Thus, a modification of the
preliminary step in the Dichotomy Algorithm enables us to
effectively solve not only a series but a single SLAE with a
Toeplitz matrix.

The structure of the article is as follows. In Section 2, we
expose the Dichotomy Algorithm for the general case. In Section 3,
we propose an economical preliminary procedure for the inversion
of tridiagonal matrices in the Toeplitz class. We consider the
problems of the stability of the dichotomy process and the
accuracy of the so-obtained solution for systems with or without
diagonal dominance. In Section 4, we give the results of numerical
experiments.  Section 5 is devoted to the summary of the work we
have done.

\section{The Dichotomy Algorithm}
Before starting the exposition of the Dichotomy Algorithm,
consider the question of mapping the data of the problem onto many
processors.
\subsection{Data decomposition}
Let  $p$ be the number of the PEs. Partition the vector of the
right-hand side and the  solution vector  ${ \bf F}$ and ${ \bf
X}$  into subvectors ${\bf Q_{i},\, U_{i}}$ as follows:

\begin{equation}
{\bf F}=\left({ \bf Q}_{1}, { \bf Q}_{2},...,{\bf Q}_{p}
\right)^\mathrm{T}=\left(f_1,f_2,...,f_{\it{size}\{\bf{
F}\}-1},f_{\it{size}\{{\bf F}\}}\right)^\mathrm{T
}\label{decom_f},
\end{equation}

\begin{equation}
{\bf X}=\left({\bf U}_{1},{\bf U}_{2},...,{\bf U}_{p}
\right)^\mathrm{T}=\left(x_1,x_2,...,x_{\it{size}\{\bf{X}\}-1},x_{\it{size}\{{\bf
X}\}}\right)^\mathrm{T}. \label{decom_x}
\end{equation}

Denoting by $\it{size}\{\bf V\}$   the number of the components of
a vector $\bf V$, require the fulfillment of the following
conditions:

$$
\begin{array}{l}
 size\{{\bf Q}_{i}\}=size\{{\bf U}_{i}\} \geq 2, \quad i=1,...,p,\\\\
\sum_{i=1}^{p}{\it size}\{{\bf Q}_{i}\}=\sum_{i=1}^{p}{\it
size}\{{\bf U}_{i}\}={\it size}\{{\bf X}\}.
\end{array}
$$

Assume that the pair of subvectors  $\left({\bf Q_{i}},{\bf
U_{i}}\right)$  belongs to the PE with the number $i$  and the row
of the matrix $A$  having number $j$  belongs to the PE containing
the pair of the elements  $\left(x_j,f_j\right)$ of
(\ref{decom_f}),(\ref{decom_x}).

Introduce in addition the following notations:
\begin{itemize}
    \item Denote the first and the last elements of a vector   $\bf V$ by  $first\{{\bf V}\}$ and $last\{{\bf V}\}$.
    \item
    Denote by  $\left\{A\right\}_{l}^{t}$ the matrix obtained from a matrix  $A$ by throwing off all rows and
    columns with the numbers less than  $l$ or greater than $t$.

    \item Denote by  $\left\{\mathbf{V}\right\}_{l}^{t}$ the subvector obtained from a vector
      $\mathbf{V}$ by throwing off the components with the numbers less than  $l$  or greater than  $t$.
    %\item Assume that
   % $\left\{A\right\}_{l}^{t}$$\left\{\mathbf{V}\right\}_{l}^{t}$=
   % $\left\{A\mathbf{V}\right\}_{l}^{t}$.
    \item  Define $\mathbf{e}_k$ as  a coordinate vector in $\mathfrak{R}^{n}$.

\end{itemize}

\subsection{The Dichotomy of a SLAE}
\begin{figure}[htb]
\center
\includegraphics[width=0.8\textwidth]{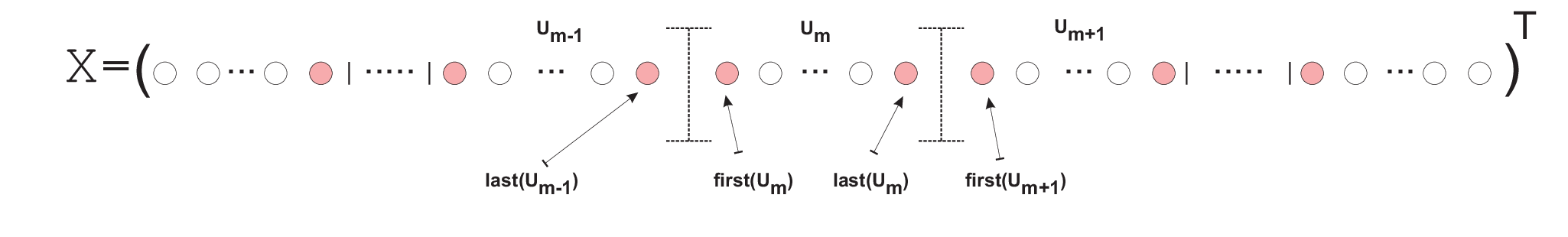}
\caption{The first and last elements of the solution vector.}
\label{pic:first_last}
\end{figure}

\begin{figure}[htb]  \center
\includegraphics[width=0.6\textwidth]{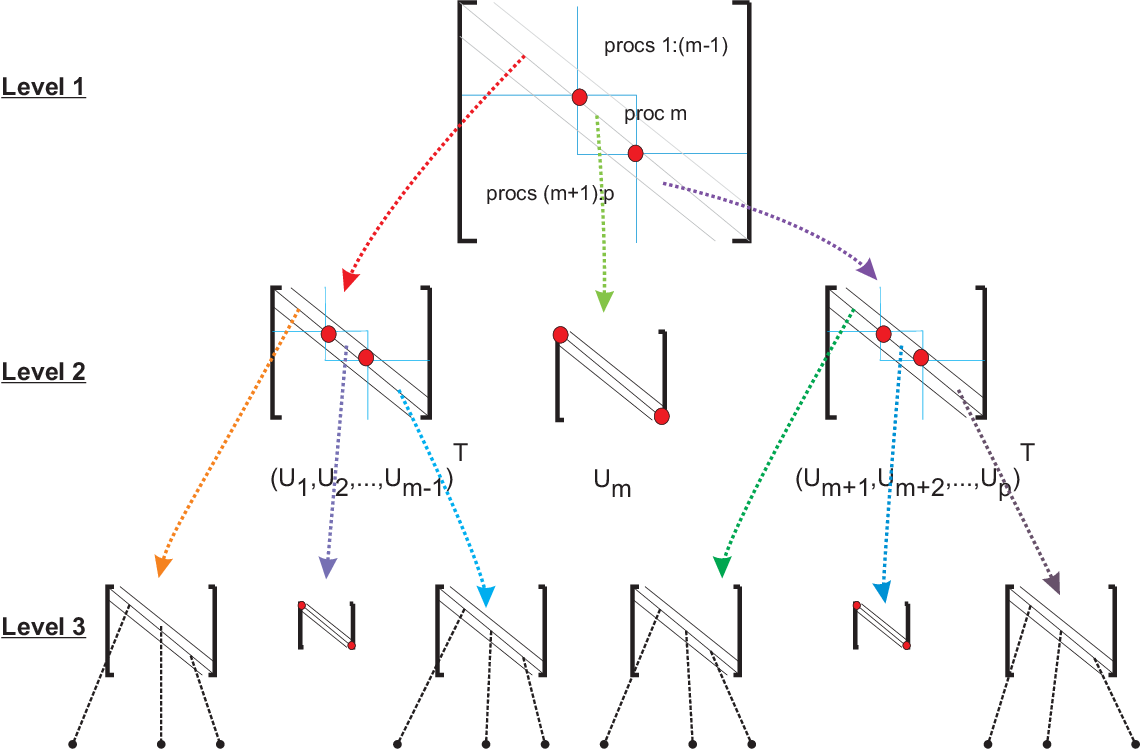}
\caption{Partition of a tridiagonal SLAE into independent
subsystems. } \label{pic:example_dichotomy}
\end{figure}

The Dichotomy Algorithm is a representative of the class of
algorithms known as "Divide \& Conquer"
\cite{Konovalov,Wang2,Lopez1994}. On each dichotomy level, the
tridiagonal system of equations obtained at the previous step is
partitioned into three independent subsystems of lesser dimensions
(Fig.~\ref{pic:example_dichotomy}) by computing the solutions in
the \mbox{$first\{{\bf U}_{m}\},\;last\{{\bf U}_{m}\}$} components
(Fig.~\ref{pic:first_last} ).

Thus, on the first dichotomy level, we construct two components of
the solution vector
\begin{equation}
\label{mrml2} \left(\mathbf{X}\right)_{m_L}\equiv
first\{\mathbf{U}_m\}, \quad \left(\mathbf{X}\right)_{m_R}\equiv
last\{\mathbf{U}_m\},
\end{equation}
situated on the  $m$th processor. The calculation of these
components makes it possible to replace the initial system by
three independent problems

\begin{equation}\label{systema}
\left\{A\mathbf{X}\right\}_{1}^{m_L-1}=\left\{\mathbf{F}\right\}_{1}^{m_L-1}-a_{m_L-1}\left(\mathbf{X}\right)_{m_L}\mathbf{e}^{\mathrm{L}},
\end{equation}

\begin{equation}\label{systemb}
\left\{A\right\}_{m_L+1}^{m_R-1} \left(
\begin{array}{c}
x_{m_L+1}\\
x_{m_L+2}\\
\dots\\
x_{m_R-1}
\end{array}\right)=
\left(
\begin{array}{c}
f_{m_L+1}-c_{m_L+1}\left(\mathbf{X}\right)_{m_L}\\
f_{m_L+2}\\
\dots\\
f_{m_R-1}-a_{m_R-1}\left(\mathbf{X}\right)_{m_R}
\end{array}\right),
\end{equation}

\begin{equation}\label{systemc}
\left\{A\mathbf{X}\right\}_{m_R+1}^{n}=\left\{\mathbf{F}\right\}_{m_R+1}^{n}-c_{m_R+1}\left(\mathbf{X}\right)_{m_R}\mathbf{e}^{\mathrm{R}}.
\end{equation}

$$
\mathbf{e}^\mathrm{R}=\left(1,0,0,...,0\right)^{\mathrm{T}},\;\mathbf{e}^\mathrm{L}=\left(0,...,0,0,1\right)^{\mathrm{T}}.
$$

On the second dichotomy level, we similarly separate systems
(\ref{systema}) and (\ref{systemc}). As a result, in   $\lceil
\log_2p \rceil$ steps, the initial system of equations
(\ref{main_eq3}) will split into $p$ independent subsystems of the
form (\ref{systemb}). The solution to (\ref{systemb}) on each
processor may be found independently with the use of any
sequential version of the Thomas algorithm
\cite{Samarski_Nikolaev,Golub1989}  for $O(size\{\mathbf{X}\}/p)$
arithmetic operations.

\subsection{The Main Formulas}
The process of the calculation of the $first\{{\bf U}_{m}\}$
,$\,last\{{\bf U}_{m}\}$  components consists in two steps: a
preliminary step carried out once for all right-hand sides of
(\ref{main_eq3}) and the dichotomy process at which the solution
is computed  for each right-hand side.

At the preliminary step, on the  $m$th processor locally without
communication interactions, we compute two rows of the matrix
$A^{-1}$

\begin{equation}
\begin{array}{cc}
\label{green_v} A^{\mathrm{T}} {\bf G}^{\mathrm{L}}_m={\bf
e}^{}_{m_L}, & A^{\mathrm{T}}{\bf G}^{\mathrm{R}}_m={\bf
e}^{}_{m_R},
\end{array}
\end{equation}

where   $m_L\,,m_R$ are defined in (\ref{mrml2}) and
$\mathbf{e}_k$ is an ort.

The vectors  $\mathbf{G}_m^{\mathrm{R,L}}$ have the sense of the
difference Green's function \cite{Samarskii2001} for the
corresponding three-point boundary value problem.

\begin{rmk}
In \cite{Terekhov:Dichotomy}, the system of linear algebraic equation with symmetric and asymmetric tridiagonal matrices were separately considered (Lemma 1 and Section 3.6). However, with allowance for the property \mbox{$(A^{\mathrm{T}})^{\mathrm{-1}}=(A^{\mathrm{-1}})^{\mathrm{T}}$}, the preliminary procedure can be essentially simplified  reducing to the calculation of the vectors $\mathbf{G}^{\mathrm{R,L}}$ following (\ref{green_v}). This makes possible, on the one hand, not to consider the case with an asymmetric matrix and, on the other hand, to exclude the situation like overflow, which may take place in explicit calculation of the inverse matrix entries.
\end{rmk}

In addition, at the preliminary step, we compute two vectors

\begin{equation}
\label{z_vector}
\begin{array}{l}
\mathbf{Z}_{m}^\mathrm{L}=\left(z^{\mathrm{L}}_{1},z^{\mathrm{L}}_{2},...,z^{\mathrm{L}}_{m_L-1},1\right)^{\mathrm{T}},\\\\
\mathbf{Z}_{m}^\mathrm{R}=\left(1,z^{\mathrm{R}}_{m_R+1},\dots,z^{\mathrm{R}}_{N-1},z^{\mathrm{R}}_{N}\right)^{\mathrm{T}},
\end{array}
\end{equation}

where the components of the vectors are found as solutions to the
systems

\begin{equation}
\left\{A\right\}_{1}^{m_L-1}\left\{\mathbf{Z}^{\mathrm{L}}
\right\}_{1}^{m_L-1}=-a_{m_L-1}\mathbf{e}^{\mathrm{L}}, \label{zl}
\end{equation}

\begin{equation}
\left\{A\right\}_{m_R+1}^{N}\left\{\mathbf{Z}^{\mathrm{R}}
\right\}_{m_R+1}^{N}=-c_{m_R+1}\mathbf{e}^{\mathrm{R}}. \label{zr}
\end{equation}

Thus, the costs of the preliminary computations for the Dichotomy
Algorithm will constitute
$size\{\mathbf{Z}_m^{\mathrm{R,L}}\}=O(N)$. Therefore, it is
appropriate to apply the Dichotomy Algorithm for solving several
SLAE with constant matrix and different right-hand sides. In this
case, the preliminary computations may be neglected.

Like in the Partition Algorithm \cite{Konovalov,Wang2}, the
Dichotomy Algorithm is based on the superposition
principle\cite{Gelphond}. That is, the components of the solution
vector are expressed via the sum of the general solution to the
homogeneous equation and a partial solution to the nonhomogeneous
equation. However, in the Dichotomy Algorithm, this principle is
realized in a somewhat different manner.

Suppose that \mbox{$\forall\,i\neq m,\;
\left\|\mathbf{Q}_i\right\|=0$}. Then the components of
(\ref{mrml2}) satisfy the identity \cite{Terekhov:Dichotomy}

\begin{equation}
\label{beta_symn}
\begin{array}{ll}
\displaystyle
\left(\mathbf{X}\right)_{m_L}=\sum_{j=m_L}^{m_R}\left(\mathbf{F}\right)_j\left(\mathbf{G}^\mathrm{L}_{m}\right)_j,&
\displaystyle
\left(\mathbf{X}\right)_{m_R}=\sum_{j=m_L}^{m_R}\left(\mathbf{F}\right)_j\left(\mathbf{G}^\mathrm{R}_{m}\right)_j.
\end{array}
\end{equation}

The remaining components of the solution vector may be determined
from the auxiliary vectors  $\mathbf{Z}_m^{\mathrm{R,L}}$  as
follows:

\begin{equation}
\left(\mathbf{X}\right)_i=\left\{
\begin{array}{ll}
\left(\mathbf{Z}_m^{\mathrm{L}}\right)_i\left(\mathbf{X}\right)_{m_L},& i \leq m_L\\\\
\left(\mathbf{Z}_m^{\mathrm{R}}\right)_i\left(\mathbf{X}\right)_{m_R},&
i \geq m_R.
\end{array}\right.
\label{solut1}
\end{equation}

If we consider all possible cases $m=1,...,p$  when
$\left\|\mathbf{Q}_m\right\|\neq 0$  while $ \forall i\neq m$
 $\left\|\mathbf{Q}_i\right\| =0$ and then sum up the
so-obtained solutions (\ref{solut1}), we come to the formula for
computing the  $first,last$-- elements

\begin{equation}
\label{theor_3} \left(\mathbf{X}\right)_k=\left\{
\begin{array}{ll}\displaystyle
\sum_{j=1}^{m-1}
\beta_{j}^\mathrm{R}\left(\mathbf{Z}_j^\mathrm{R}\right)_k+\sum_{j=m}^{p}
\beta^\mathrm{L}_{j}\left(\mathbf{Z}_j^\mathrm{L}\right)_k, &
\left(\mathbf{X}\right)_k=first\{\mathbf{U}_m\},\\
\\\\
\displaystyle   \sum_{j=1}^{m}
\beta_{j}^\mathrm{R}\left(\mathbf{Z}_j^\mathrm{R}\right)_k+\sum_{j=m+1}^{p}
\beta^\mathrm{L}_{j}\left(\mathbf{Z}_j^\mathrm{L}\right)_k, &
\left(\mathbf{X}\right)_k=last\{\mathbf{U}_m\},
\end{array}\right.
\end{equation}

where

\begin{equation}
\label{beta_symn}
\begin{array}{ll}
\displaystyle
\beta_m^\mathrm{L}=\sum_{j=m_L}^{m_R}\left(\mathbf{F}\right)_j\left(\mathbf{G}^\mathrm{L}_{m}\right)_j,&
\displaystyle
\beta_m^\mathrm{R}=\sum_{j=m_L}^{m_R}\left(\mathbf{F}\right)_j\left(\mathbf{G}^\mathrm{R}_{m}\right)_j.
\end{array}
\end{equation}
Here the indices $m_R,m_L$  are found locally on each processor by
(\ref{mrml2}).

The vectors
$u_m,v_m$  are the solution to the interior difference boundary
value problem with respect to the subvector $\mathbf{U}_m$,
whereas the vectors $\mathbf{Z}_m^{\mathrm{R,L}}$ are the solution
to the exterior boundary value problem. This difference is of
principle. So, in the realization of the Dichotomy Algorithm, the
computation of the $first,last$ -- components is reduced to the
calculation of the sums (\ref{theor_3}), whereas, in the Partition
Algorithm, it is necessary to solve a "reduced"\
\cite{Konovalov,Wang2} SLAE of dimension $2p-2$ by means of some
parallel variant of the Gauss Elimination Method. Since the
computation of the sums on a multiprocessor computing system can
be realized with a greater efficiency then the Gauss elimination
method, the Dichotomy Algorithm makes it possible to reach a
greater productivity for the problems of the form (\ref{main_eq3})
than the Partition Algorithm.

\subsection{MPI realization of the Dichotomy Algorithm}
\begin{figure}[htb]
 \center
\includegraphics[width=0.65\textwidth]{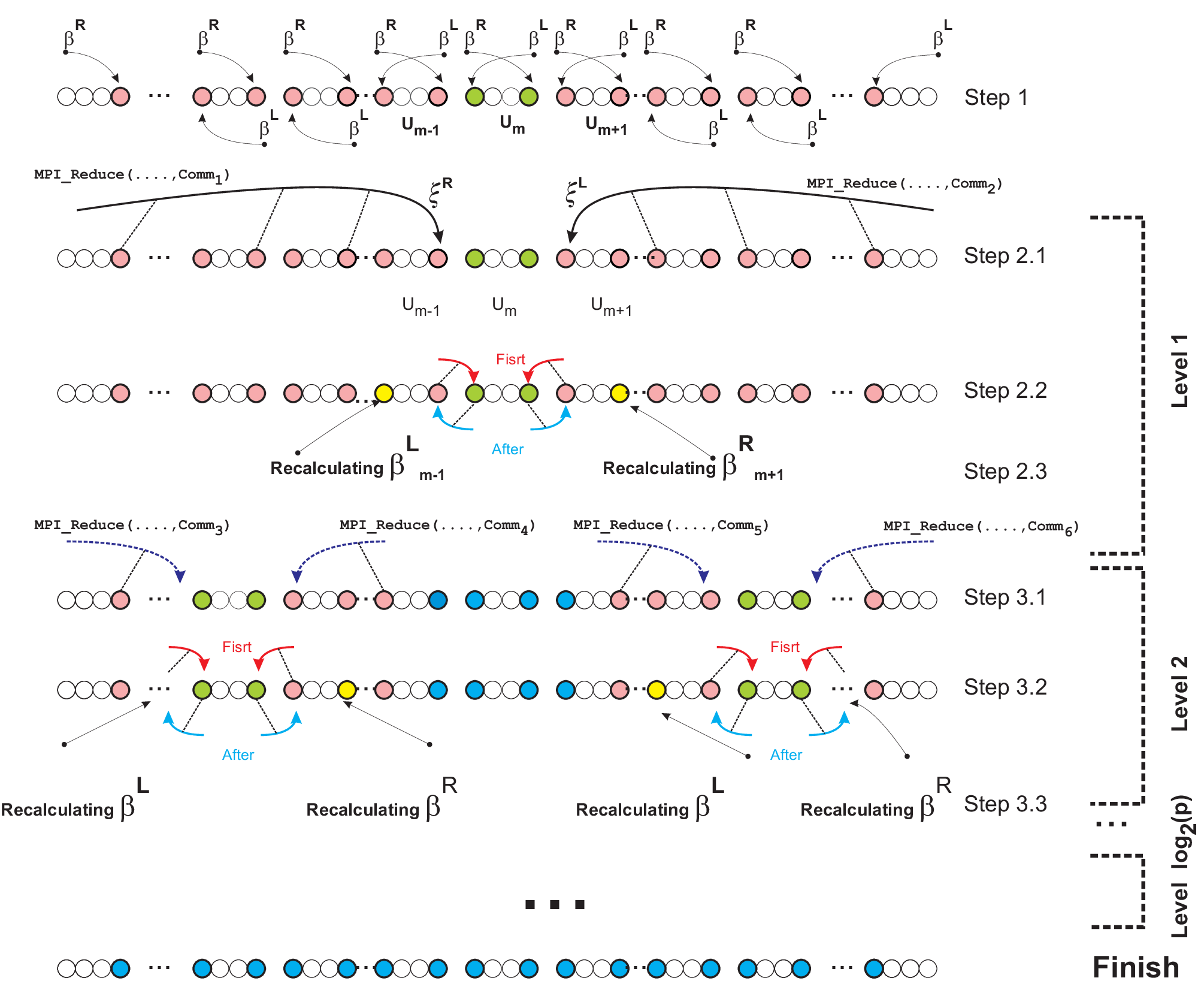} \caption{The MPI-realization of the Dichotomy Process.} \label{pic_analitical}
\end{figure}

 Fig.~\ref{pic_analitical} contains the MPI-realization of (\ref{theor_3}) for partitioning the system
\mbox{$A\mathbf{X}=\mathbf{F},$} \mbox{$size\{\mathbf{F}\}=N$}.

Assume that the vectors
$\mathbf{G}_m^{\mathrm{R,L}},\mathbf{Z}_m^\mathrm{R,L}$ are
already defined at the preliminary step. Then, at step 1 of the
Dichotomy Process, by (\ref{beta_symn}) each processor computes
two local quantities $\beta_m^{\mathrm{R,L}}$, where $m$ is the
processor number. The arithmetic costs at this step are
$O(size\{\mathbf{U}_i\})$ for each PE. In the sequel, in the
separation of the systems, the quantities $\beta_m^{\mathrm{R,L}}$
will be recomputed for a number of operations of order $O(1)$.

At the first dichotomy level, consisting of steps 2.1--2.2--2.3,
we compute two components of (\ref{mrml2}) with the use of
(\ref{theor_3}) as follows.

At step 2.1, by calling the collective function $mpi\_Reduce$ over
the communicators $\mathrm{Comm_1},\mathrm{Comm_2}$,  we compute
the quantities

\begin{equation}
\label{step21}
\begin{array}{ll}\displaystyle
\xi^{\mathrm{R}}=\sum_{j=1}^{m-1}
\beta_{j}^\mathrm{R}\left(\mathbf{Z}_j^\mathrm{R}\right)_{m_L}, &
\displaystyle \xi^{\mathrm{L}}=\sum_{j=m+1}^{p}
\beta^\mathrm{L}_{j}\left(\mathbf{Z}_j^\mathrm{L}\right)_{m_R}.
\end{array}
\end{equation}

At step 2.2, processor  $m-1$ sends\footnote{This  exchange is
designated as "First"\ in  Fig.~\ref{pic:example_dichotomy}.}  to
processor $m$  the quantity $\xi^{\mathrm{R}}$ and processor
$m+1$, the quantity $\xi^{\mathrm{L}}.$

Now, the sought components situated on the  $m$th processor may be
computed as \cite{Terekhov:Dichotomy}

\begin{equation}
\label{theor33}
\begin{array}{ll} \displaystyle
\left(\mathbf{X}\right)_{m_L}=first\{\mathbf{U}_m\}=\xi^{\mathrm{R}}+\xi^{\mathrm{L}}\frac{\left(\mathbf{G}^{\mathrm{L}}_m\right)_{m_R}}{\left(\mathbf{G}^{\mathrm{R}}_m\right)_{m_R}}+\beta_{m}^{\mathrm{L}},\\\\
\displaystyle
\left(\mathbf{X}\right)_{m_R}=last\{\mathbf{U}_m\}=\xi^{\mathrm{R}}\frac{\left(\mathbf{G}^{\mathrm{R}}_m\right)_{m_L}}{\left(\mathbf{G}^{\mathrm{L}}_m\right)_{m_L}}+\xi^{\mathrm{L}}+\beta_{m}^{\mathrm{R}}.
\end{array}
\label{step23}
\end{equation}

\begin{rmk}
The quantities
${\left(\mathbf{G}^{\mathrm{R}}_m\right)_{m_L}}/{\left(\mathbf{G}^{\mathrm{L}}_m\right)_{m_L}}$
and
${\left(\mathbf{G}^{\mathrm{L}}_m\right)_{m_R}}/{\left(\mathbf{G}^{\mathrm{R}}_m\right)_{m_R}}$
make it possible to "transfer"\ the boundary condition from the
component  $first$  to the component $last$  and vice versa, i.e.,
inside the vector $\mathbf{U}_m$.
\end{rmk}

Furthermore, to exclude the thus-found components from the system
of equations,  processor $m$  sends\footnote{This exchange is
designated as "After" in Fig.~\ref{pic:example_dichotomy}.}  to
processor $m-1$  the quantity $\delta^{\mathrm{L}}=\left(\xi^{\mathrm{L}}\frac{\left(\mathbf{G}^{\mathrm{L}}_m\right)_{m_R}}{\left(\mathbf{G}^{\mathrm{R}}_m\right)_{m_R}}+\beta_{m}^{\mathrm{L}}\right)\left(\mathbf{Z}_{m}^\mathrm{L}\right)_{m_L-1}$ and to processor $m+1$, the quantity

$\delta^{\mathrm{R}}=\left(\xi^{\mathrm{R}}\frac{\left(\mathbf{G}^{\mathrm{R}}_m\right)_{m_L}}{\left(\mathbf{G}^{\mathrm{L}}_m\right)_{m_L}}+\beta_{m}^{\mathrm{R}}\right)\left(\mathbf{Z}_{m}^\mathrm{R}\right)_{m_R+1}$
respectively.

At step 2.3, the vector of the right-hand side of SLAE
(\ref{systema}),(\ref{systemb}),(\ref{systemc}) is modified, and
the processors with numbers \mbox{$m-1,m+1$}  recompute the
quantities

\begin{equation}
\label{recalc_beta}
\begin{array}{cc}
\hat{\beta}_{m-1}^\mathrm{R}=\beta_{m-1}^\mathrm{R}+\delta^{\mathrm{L}},&
\hat{\beta}_{m+1}^\mathrm{L}=\beta_{m+1}^\mathrm{L}+\delta^{\mathrm{R}},\\\\
\hat{\beta}_{m-1}^\mathrm{L}=\beta_{m-1}^\mathrm{L}+\delta^{\mathrm{L}}\frac{\left(\mathbf{G}^{\mathrm{L}}_{m-1}\right)_{t_R}}{\left(\mathbf{G}^{\mathrm{R}}_{m-1}\right)_{t_R}},&
\hat{\beta}_{m+1}^\mathrm{R}=\beta_{m+1}^\mathrm{R}+\delta^{\mathrm{R}}\frac{\left(\mathbf{G}^{\mathrm{R}}_{m+1}\right)_{q_L}}{\left(\mathbf{G}^{\mathrm{L}}_{m+1}\right)_{q_L}},\\\\

\end{array}
\end{equation}
where $\left(\mathbf{X}\right)_{t_R}\equiv
last\{\mathbf{U}_{m-1}\}, \quad \left(\mathbf{X}\right)_{q_L}\equiv
first\{\mathbf{U}_{m+1}\},
$

On the next dichotomy level, an analogous separation process is
applied to subsystems (\ref{systema}),(\ref{systemc}), where the
quantities (\ref{recalc_beta}) are used instead of
$\beta_{m-1,m+1}^{\mathrm{R,L}}$. Moreover, the steps having
numbers $s.1,s.2,s.3$, (where $s=1,2,...,\lceil \log_2p \rceil$ )
are equivalent to steps $2.1,2.2,2.3$. However, on the $s$th
level, $2^{s-1}$ already independent systems obtained at the
previous step are separated independently.

\subsection{Optimization of the Communication Costs}
\label{Par:cost} Since the number of arithmetic operations at the
step of solving the "reduced" SLAE is relatively small, the
communication cost of the algorithm of solving this subproblem
defines the efficiency of the Partition Algorithm. The "reduced"\
system may be solved by the Cyclic Reduction Method
\cite{Tridiagonal5,Johnsson:cyclic}. In \cite{Terekhov:Dichotomy}
it was shown that the estimates of the computation time for the
Dichotomy Algorithm and the Cyclic Reduction coincide in order if,
first, the delay time $\alpha$ before the data transfer is
insignificant and, second, not one but several series of problems
are solved simultaneously

$$
\label{dich_time2}
\begin{array}{l}
{T}^{Dichotomy}_{p}
=\alpha\left[\log_2(p)+1\right]\log_2(p)+l\left(\log_{2}(p)-\frac{p-1}{p}\right)\left(\gamma+2\beta\right)\approx\\\\
\approx
2\log_{2}(p)\left(l\beta+l\gamma/2\right),\\\\
{T}^{Cyclic\,Reduction}_{p}=2\log_2(p)\left(\alpha+l\beta+l\gamma\right)\approx
2\log_2(p)\left(l\beta+l\gamma\right).
\end{array}
$$

Here $\beta$  is the time of transfer of a real number from one PE
to another,  $\gamma$ is the time spent on the operation of
addition of two numbers, $l$  is the number of the series of the
simultaneously solved SLAEs.

In comparison with the Partition Algorithm $+$ the Cyclic
Reduction, the Dichotomy Algorithm gives a higher speedup
coefficient due to lesser time costs for the synchronization of
the computations as well as on the data exchange between the PEs.

The reduction of the communication time is possible due to the
fact that the main communication operation ("+") of the Dichotomy
Algorithm is associative. The architecture of modern
supercomputers is such that the time of pairwise interactions may
differ substantially for different processors
\cite{Parhami,Dongarra:Reduce,MPI2}. The associativity of the
computations makes it possible to define the order of the
interactions of the PEs  on the level of a communication  library
or a programming language so as to minimize the time of data
exchanges by  taking into account the architecture of the
supercomputer.

The reduction of the interprocessor exchange time is not a
sufficient condition for the efficiency of the parallel algorithm
since it is necessary to take into account the presence of
synchronization of the computations. It often happens that most of
the PEs expect the computation results from several PEs, which
leads to the decrease of the total efficiency due to down-time
periods of the computational resources.   In the context of the
Dichotomy Algorithm, the problem of minimization of the
synchronization time is solved as follows. First, the number of
groups of interacting processors increases with the number of the
dichotomy level but the number of processes in each of them
decreases; hence, so does the time of interprocessor exchanges.
Second, the Dichotomy Algorithm almost does not require to
synchronize the computations in passing from one dichotomy level
to another (especially, on the first levels, where the number of
communications is maximal). Thus, for beginning the process of
partitioning the system of equations on dichotomy level  $s+1$, it
is first necessary to modify the quantities (\ref{recalc_beta}) at
step $s.3$. However, the processors that do not compute the
 $\hat{\beta}_{m}^{\mathrm{R,L}}$ at step  $s.3$  may start the summation of the series (\ref{step21}) at step $(s+1).1$
without waiting for the result of the previous steps. Thus, on a
majority of the processors, steps $s.2,s.3$, and $(s+1).1$  may be
carried out with a high degree of parallelism.

\section{The Parallel Dichotomy Algorithm for Toeplitz Matrices}
\label{section_toeplitz} The realization of the Dichotomy
Algorithm on a parallel computational system may be regarded as
the solution of two separate subproblems:

\begin{enumerate}
    \item Minimization of the number of arithmetic operations in the computation of the auxiliary vectors  $\mathbf{Z}_m^{\mathrm{R,L}},\,
\mathbf{G}_m^{\mathrm{R,L}}$.
    \item Minimization of the time of the communications and synchronization in the realization of
    (\ref{theor_3}).
\end{enumerate}

Since, at the second step of the Dichotomy Algorithm, the form of
the SLAE does not matter, we, in order to guarantee the
possibility of solving not only a series but a single equation
with matrix of the form (\ref{Toelitz_matrix}), pose the problem
of reducing the number of arithmetic operations at the preliminary
step. Note that, in view of (\ref{theor_3}),(\ref{beta_symn}) we
do not need to determine all components of the auxiliary vectors
but only those used. Thus, we need to find only $O(N/p)$
components of $\mathbf{G}_m^{\mathrm{R,L}}$ and $\lceil \log_2
p\rceil$ components of ${\bf Z}_m^\mathrm{R,L}$.

\subsection{Optimization of the Preliminary Computations}
For a general matrix  (\ref{main_eq3}), finding the necessary
components $z_i^{\mathrm{R,L}},g_i^{\mathrm{R,L}}$ from
(\ref{green_v}),(\ref{z_vector}) requires $O(N)$ arithmetic
operations. For Toeplitz matrices, the components of the auxiliary
vectors may be found in a lesser number of operations in
accordance with the following theorem.

\begin{thm}\cite{Samarski_Nikolaev,MatrixInvert,Gelphond} Assume that we need to solve a SLAE $$T {\bf Y}= { \bf F}$$ with
matrix (\ref{Toelitz_matrix}) of order $N$. Then the  $n$th
component of the solution may be calculated as
\begin{equation}
\label{theor1}
y_n=\sum_{k=1}^{n-1}\frac{(q_1q_2)^{n-k}\left(q_2^{N+1-n}-q_1^{N+1-n}\right)\left(q_2^k-q_1^k\right)}{(q_2-q_1)\left(q_2^{N+1}-q_1^{N+1}\right)}\cdot
\frac{f_k}{t_1}+\sum_{k=n}^{N}\frac{\left(q_2^{N+1-k}-q_1^{N+1-k}\right)\left(q_2^n-q_1^n\right)}{(q_2-q_1)\left(q_2^{N+1}-q_1^{N+1}\right)}\cdot
\frac{f_k}{t_1},
\end{equation}

$$q_1=\frac{-t_0+\sqrt{t_0^2-4t_{-1}t_1}}{2t_1},\quad
q_2=\frac{-t_0-\sqrt{t_0^2-4t_{-1}t_1}}{2t_1}.$$ Moreover, the
solution does not exist if $q_1^{N+1}=q_2^{N+1}$ but $q_1\neq
q_2$.
\end{thm}

Involving the fact that the right-hand sides in
(\ref{green_v}),(\ref{zl}),(\ref{zr}) may contain a unique nonzero
component, it is easy to prove that, for computing one component
of the vectors
$\mathbf{Z}_m^{\mathrm{R,L}},\mathbf{G}_m^{\mathrm{R,L}}$, by
(\ref{theor1}) we will need $O(1)$ arithmetic operations.
Therefore, at the preliminary step as well as at the step of the
partition of the SLAE, it will be necessary to carry out
$O(N/p+\log_2p)$ arithmetic operations. Thus, the Dichotomy
Algorithm for Toeplitz matrices may be applied not only to a
series but to a single problem.

A similar economic preliminary procedure may be constructed for
matrices of a somewhat more general kind than
(\ref{Toelitz_matrix}). For example, in \cite{MatrixInvert,Meurant1992}, explicit
expressions are given for computing the elements of the inverse
matrix to a matrix of the form
\begin{equation}
\label{Toelitz_matrix2}
 A=\left(
\begin{array}{ccccc}
  t_0+\psi & t_1  &  &  &  \Large 0\\
  t_{-1} & t_0 & t_1  &  &  \\
%   & t_{-1} & t_0 & t_1 &  &  \\
    & \ddots & \ddots & \ddots &  \\
    &  & t_{-1} & t_{0} & t_{1} \\
   0&  &  & t_{-1} & t_0+\chi \\
\end{array}
\right).
\end{equation}

Thus, if, for some tridiagonal matrix, it is possible to compute a
separate component of the vectors
$\mathbf{Z}_m^{\mathrm{R,L}},\mathbf{G}_m^{\mathrm{R,L}}$ for
$O(1)$ arithmetic operations then, in this case,  it is possible
to effectively solve a series as well as a single problem.

\subsection{Inversion of the Operator  $\nabla^2$}
In solving the first boundary value problem for the Poisson
equation by such methods as the Multigrid Line Relaxation, the
Variable-Separation Method etc., it is necessary to solve a series
of equations of the form

\begin{equation}
\label{razdiffeq} \left\{
\begin{array}{ll}\displaystyle
\frac{y_{i+1}-2y_{i}+y_{i-1}}{h^2}+\lambda_k y_i=-f_i,\;
\lambda_k\in R & 1\leq i\leq N-1,\quad
k=1,...,M,\\\\y_0=\mu_1,\quad y_N=\mu_2.
\end{array}\right.
\end{equation}
Consider an economic algorithm for computing of the components of
the vectors
$\mathbf{Z}_m^\mathrm{R,L},\mathbf{G}_m^{\mathrm{R,L}}$  for
problems (\ref{razdiffeq}). It is known
\cite{Samarski_Nikolaev,MatrixInvert,Fonseca2001} that the
solution to (\ref{razdiffeq}) is given by (\ref{theor1}) and may
be written down in the form
\begin{equation}
\label{chebU}
y_n=\frac{U_{N-n-1}(x)}{U_{N-1}(x)}\left[\mu_1+\sum_{k=1}^{n-1}U_{k-1}(x)f(k)\right]+
\frac{U_{n-1}(x)}{U_{N-1}(x)}\left[\mu_2+\sum_{k=n}^{N-1}U_{N-k-1}(x)f(k)\right],
\end{equation}
\\\\
where $x=1-h^2\lambda/2\neq \cos{\frac{k\pi}{N}},\; k=1,2,...,N-1$
and  $U_n(x)$ is the Chebyshev polynomial of the second kind of
$n$ degree \cite{abramowitz+stegun,Mason:Cheb}.

With (\ref{mrml2}) taken into account, the solution to
(\ref{z_vector}) on the  $m$th processor has the form

\begin{equation}
\label{Zdiff}
\begin{array}{ll}
z^{\mathrm{R}}_{i}=\frac{U_{N-i-1}(x)}{U_{N-m_R-1}(x)},& m_R\leq i\leq N,\\\\
z^\mathrm{L}_{i}=\frac{U_{i-1}(x)}{U_{m_L-1}(x)},&  1\leq i\leq
m_L.
\end{array}
\end{equation}

The solution to (\ref{green_v}) has the form
\begin{equation}
\label{Greendiff}
\begin{array}{ll}
g^{\mathrm{L}}_{i}=
\frac{U_{N-i-1}(x)}{U_{N-1}(x)}U_{m_L-1}(x),& m_L\leq i\leq m_R,\\\\
g^{\mathrm{R}}_{i}=\frac{U_{i-1}(x)}{U_{N-1}(x)}U_{N-m_R-1}(x),&
m_L\leq i\leq m_R.
\end{array}
\end{equation}

Since the necessary components of  $\mathbf{G}_m^{\mathrm{R,L}}$
are situated  successively (\ref{beta_symn}), having computed
$g^{\mathrm{R,L}}_{m_R,m_L}$ by (\ref{Greendiff}), it is more
economic to compute the remaining components by solving the
systems

\begin{equation}\label{system_green1}
\left\{A\mathbf{G}^\mathrm{R}\right\}_{m_L+1}^{m_R-1}=-t_{-1}g^\mathrm{R}_{m_L}\mathbf{e}^{\mathrm{R}}-t_1g^\mathrm{R}_{m_R}\mathbf{e}^{\mathrm{L}},
\end{equation}

\begin{equation}\label{system_green2}
\left\{A\mathbf{G}^\mathrm{L}\right\}_{m_L+1}^{m_R-1}=-t_{-1}g^\mathrm{L}_{m_L}\mathbf{e}^{\mathrm{R}}-t_1g^\mathrm{L}_{m_R}\mathbf{e}^{\mathrm{L}}.
\end{equation}

\subsection{Accuracy analysis}
Since the values of  $U_k(x)$  outside the interval $[-1,1]$ begin
to increase rapidly with  $k$, we cannot exclude an overflow-type
situation in a program realization of
(\ref{Zdiff}),(\ref{Greendiff}) for  $|x|>1$
\cite{IEEE_FLOAT,Higham:accuracy}. To overcome this difficulty, we
do the following. Let $N_0$ be the degree of the polynomial
greater than which the quantity $U_k(x),\; k>N_0,\, |x|>1$ cannot
be computed because the result exceeds the boundary of the
admissible values of the real variable. A Chebyshev polynomial of
the second kind admits the representation \cite{abramowitz+stegun}
\begin{equation}
\label{chebysh}
U_n(x)=\frac{1}{2\sqrt{x^2-1}}\left[\left(x+\sqrt{x^2-1}\right)^{n+1}-\left(x+\sqrt{x^2-1}\right)^{-(n+1)}\right],
\quad |x|\geq 1.
\end{equation}

Then the relation
\begin{equation}
\label{ineq1}
\left(x+\sqrt{x^2-1}\right)^{-(n+1)}\ll\left(x+\sqrt{x^2-1}\right)^{(n+1)},
\end{equation}
holds  $\forall \, n>N_0$, from which $\forall \, n>N_0$   we may
assume with a good accuracy  that
\begin{equation}
\label{chebysh2}
U_n(x)\simeq\frac{1}{2\sqrt{x^2-1}}\left(x+\sqrt{x^2-1}\right)^{n+1}.
\end{equation}

Now, substituting (30) into (25) and collecting similar summands,
we infer
\begin{equation}
\label{Greendiff3} \displaystyle \begin{array}{ll}
g^{\mathrm{L}}_i\simeq\frac{1}{2\sqrt{x^2-1}}\left[\eta^{(m_L-i)}+\eta^{(-2N-i+m_L)}-\eta^{(-2N+m_L+i)}-\eta^{(-m_L-i)}\right],&m_L\leq i,\\\\
g^{\mathrm{R}}_i\simeq\frac{1}{2\sqrt{x^2-1}}\left[\eta^{(i-m_R)}+\eta^{(-2N+i-m_R)}-\eta^{(-2N+m_R+i)}-\eta^{(-i-m_R)}\right],&
m_R \geq i,
\end{array}
\end{equation}
where $\eta=x+\sqrt{x^2-1}$.

Since the exponents in (\ref{Greendiff3}) are at most zero and
$|x|>1$, the situation of an overflow of variables is excluded.

\begin{rmk}
In solving tridiagonal systems of dimension less than $N_0$, the
computations of the components of the auxiliary vectors should be
performed by (\ref{Zdiff}),(\ref{Greendiff}), (\ref{chebysh}) so
as to avoid loss of accuracy because of a possible violation of
(\ref{ineq1}).
\end{rmk}

\begin{rmk}
If  $|x|<1$ then $U_k(x)\leq k+1$. Hence no overflow of variables
arises for $x\neq \cos{\frac{k\pi}{N}},\; k=1,2,...,N-1$  and
reasonable $N$.
\end{rmk}

\subsection{Stability analysis}
It is proved in \cite{Terekhov:Dichotomy} that a sufficient
condition of stability for the Dichotomy Algorithm is given by the
diagonal dominance of the matrix of the SLAE

\begin{equation}
\label{diag1} \left|b_i\right| \geq
\left|a_i\right|+\left|c_i\right|,\quad i=2,...,N-1,
\end{equation}

\begin{equation}
\label{diag2} \left|b_1\right| \geq \left|a_1\right|,\quad
\left|b_N\right| \geq \left|c_N\right|,
\end{equation}

and at least one of the inequalities is strict.

\begin{figure}[!h]
\includegraphics[width=0.5\textwidth]{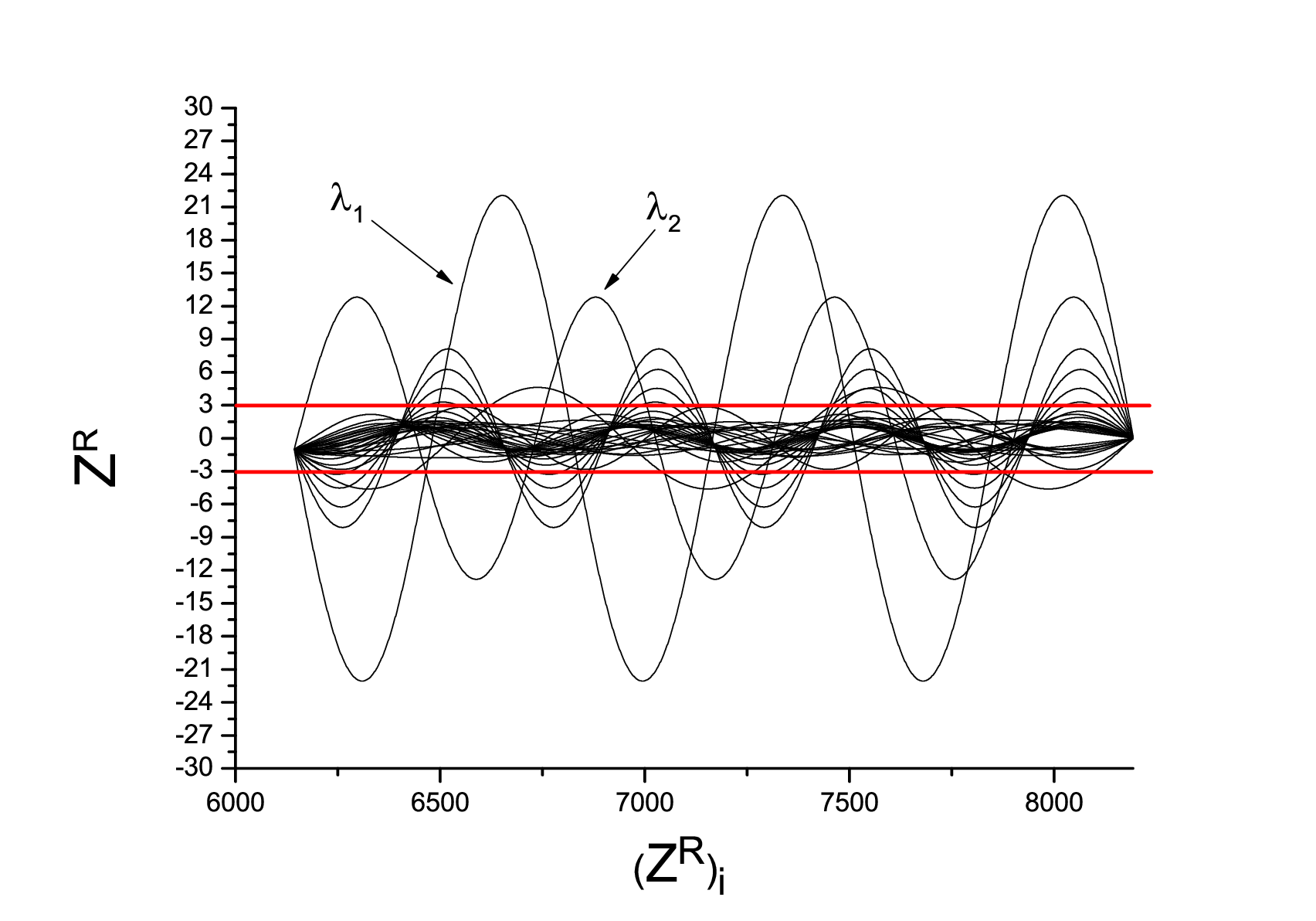} \hfill
\includegraphics[width=0.5\textwidth]{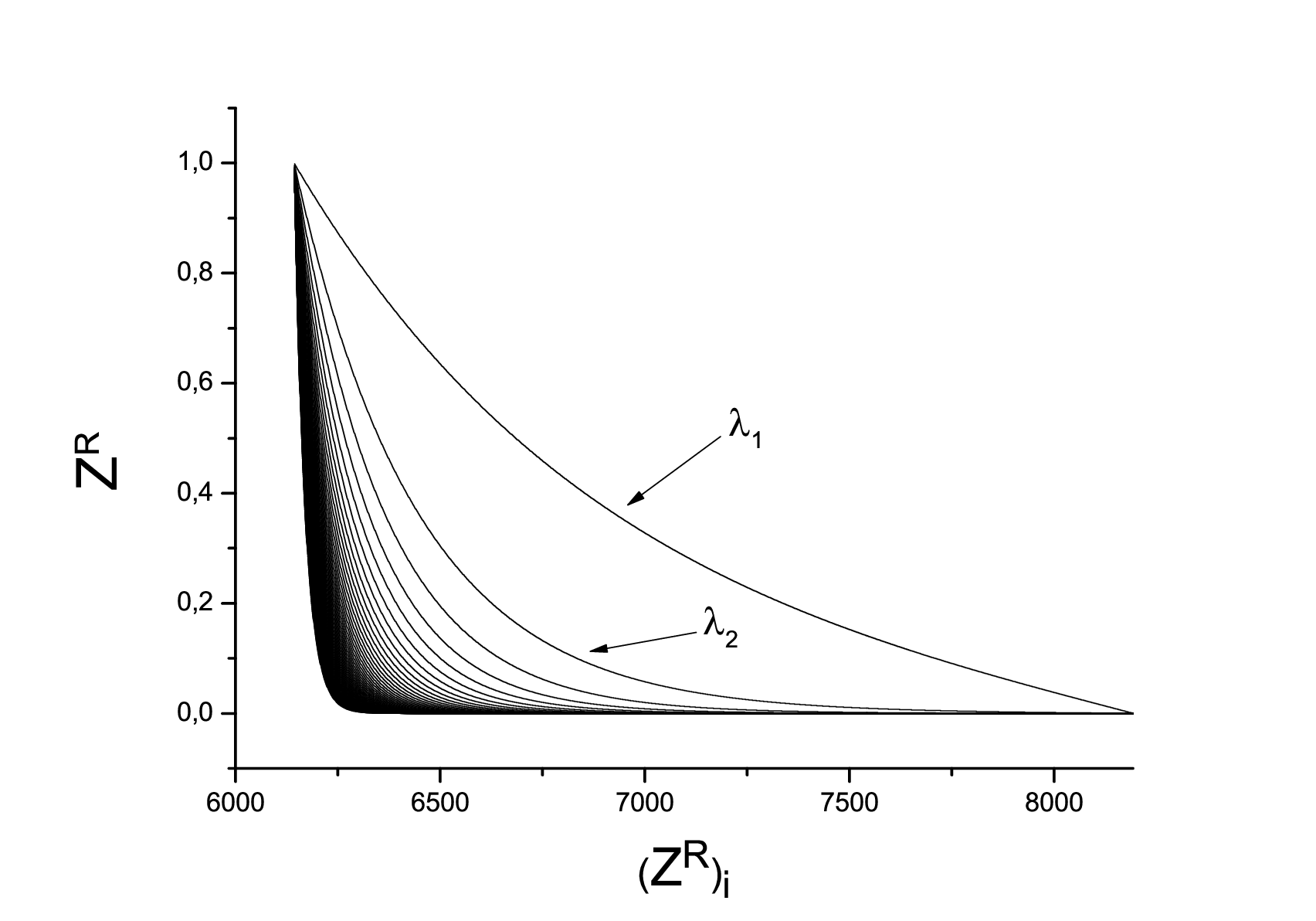}\\
\parbox{0.5\textwidth}{\center a) $\lambda_k >\lambda_{k+1}> 0$}\hfill
\parbox{0.5\textwidth}{\center b) $\lambda_{k+1} <\lambda_k < 0$}
\caption{ The values of the components of
$\mathbf{Z}_2^\mathrm{R}$ for the case  $N=8192,\,p=4$ in solving
(\ref{razdiffeq}) for different values $\lambda$.}
\label{pic:ZRZL}
\end{figure}

For problem (\ref{razdiffeq}) in the case of $\lambda \leq 0$, the
matrix of the SLAE has diagonal dominance, which guarantees the
stability of the Dichotomy Algorithm. However, for $\lambda>0$,
there is no diagonal dominance and accumulation of roundoff error
can happen. This follows from the fact that, by the failure of the
maximum principle \cite{Samarskii2001,Strikwerda2004}, the
components of the vectors $\mathbf{Z}_m^{\mathrm{R,L}}$ may have
their modulus greater than one (Fig.~\ref{pic:ZRZL}). In this
case, the error of the computation of the $first,last$-components
committed on the  $s$th dichotomy level passes to  $s+1$ level by
means of (\ref{recalc_beta}) and is then "strengthened"  by
multiplying by $\left|\mathbf{Z}_m^{\mathrm{R,L}}\right|\leq
\gamma$ in (\ref{theor_3}). Thus, if the number of the dichotomy
levels is great and $\gamma$ is much greater than 1 then an
accuracy loss is possible.

If we take the estimate
\begin{equation}
\varepsilon_i=\prod_{j=1}^{l}\left
|\left(\mathbf{Z}_{n_j}^{\mathrm{R,
L}}\right)_i\right|\varepsilon\leq\gamma^l\varepsilon, \label{err}
\end{equation}
as the maximal error of the computation of the component
$\left(\mathbf{X}\right)_i$, where $l$ is the number of the
dichotomy step at which the sought component will be computed,
$\varepsilon$ is the initial error. Then, for example, for
$p=4096,l=12$  and $\gamma \leq 3$, the computation error at worst
goes six digits to the right, which is acceptable for many
problems. Thus, for $\gamma>1$, we need to perform an a priori
(\ref{err}) as well as an a posteriori check of the accuracy of
the so-obtained solution.
\subsection{Inversion of the Operator $\nabla^2$ for 3D Case.}
For solving the 3D problems, it is necessary to provide the possibility of using a large number of processors for algorithms efficient with respect to the number of arithmetic actions. From the viewpoint of data decomposition between the processors, the best choice is the variable directions method (ADI) \cite{Samarski_Nikolaev}, where at each step of the iteration process the tridiagonal SLAEs should be solved for all the three directions, which allows one to use a large number of processors within the limits of one calculation. Unfortunately, compared to the 2D case, the 3D-ADI algorithm is not efficient because of its low convergence rate \cite{Samarskii2001}. Thus, via the Fourier transform it is reasonable to solve not a 3D problem but a series of independent 2D problems:
\begin{equation}
\label{razdiffeq22}
\begin{array}{lr}\displaystyle
%\frac{(l)_{i+1,j}-2y(l)_{i,j}+y(l)_{i-1,j}}{h_x^2}+\frac{y(l)_{i,j+1}-2y(l)_{i,j}+y(l)_{i,j-1}}{h_y^2}
\frac{y(l)_{i+1,j}-2y(l)_{i,j}+y(l)_{i-1,j}}{h_1^2}+\frac{y(l)_{i,j+1}-2y(l)_{i,j}+y(l)_{i,j-1}}{h_2^2}-\frac{4}{h_3^2}\sin^2\left(\frac{\pi(l-1)}{2(N_3-1)}\right)y(l)_{i,j}=-f(l)_{i,j},\\\\
i=2,...,N_1-1,\; j=2,...,N_2-1,\;l=2,...,N_3-1. %,\\\\y_0=\mu_1,\quad y_N=\mu_2.
\end{array}
\end{equation}
For solving the problem (\ref{razdiffeq22}), we can use the parallel version of the variable separation method \cite{Terekhov:Dichotomy}. However, this algorithm drastically limits the number of used processors because the tridiagonal SLAEs are solved only in one direction and, consequently, only the one-dimensional data decomposition will be efficient because for other directions it is necessary to use the fast Fourier transform algorithm. In order to overcome the problems of data decomposition and ensure efficiency of the method in terms of the number of operations, for solving the problems (\ref{razdiffeq22}), one should use the 2D-ADI method \cite{Terekhov:Dichotomy}. The calculation efficiency is achieved due to the fact that as the value of $l$ increases from $2$ to $N_3-1$ the condition number for an equation of the form (\ref{razdiffeq22}) decreases as well as the necessary number of iterations of the ADI algorithm, Fig. \ref{pic:3d}.a. Thus, the hybrid parallel algorithm under consideration is comparable in efficiency with the variable separation method, however, contrary to the latter, it allows 2D efficient data decomposition between the processors.
\section{Numerical Experiments}
For estimating the efficiency of the above-proposed modification
of the dichotomy algorithm for Toeplitz matrices, consider the
following boundary value problem
\begin{equation}
\label{poisson_test}
\begin{array}{lr}
\triangle u=f(x) ,\quad x \in G \subset \mathfrak{R}^\alpha,\;\alpha=2,3 ,& \left. u\right | _\Gamma=0.
\end{array}
\end{equation}
For solving the difference problem in $Fortran-90$  with the use
of the communication MPI--library, we have realized the variable
separation method following the scheme of
\cite{Terekhov:Dichotomy}.  For the 2D problem the approximation of
(\ref{poisson_test}) was carried out with the second order of
accuracy on a uniform mesh with the number of the points
\mbox{$N_1=N_2=2^k$}, \mbox{$k=13...16$} and with the number of the points \mbox{$N_1=N_2=N_3=2^k$}, \mbox{$k=9...13$} for the 3D problem. For estimating the cost of the
preliminary computations, the calculation of the vectors
$\mathbf{G}^{\mathrm{R,L}}_m,\, \mathbf{Z}^{\mathrm{R,L}}_m$ for
problem (\ref{poisson_test}) was realized in two variants: by
(\ref{green_v}),(\ref{z_vector}) for an arbitrary
symmetric matrix and by  (\ref{Greendiff3}),(\ref{z_vector})  for Toeplitz
matrices. Since the matrix of the SLAE is not used at the second
stage of the dichotomy algorithm, one MPI-realization of
(\ref{theor_3}),(\ref{beta_symn}) is enough for these cases. The calculation was performed on the "Lomonosov"\ supercomputer of Moscow State University. The supercomputer comprises Intel Xeon X5570 four-core processors operating at 2.93 GHz in the Infiniband QDR communication environment. Each computational node contains two processors and $12$ GB of RAM.

\begin{figure}[!h]
\includegraphics[width=0.5\textwidth]{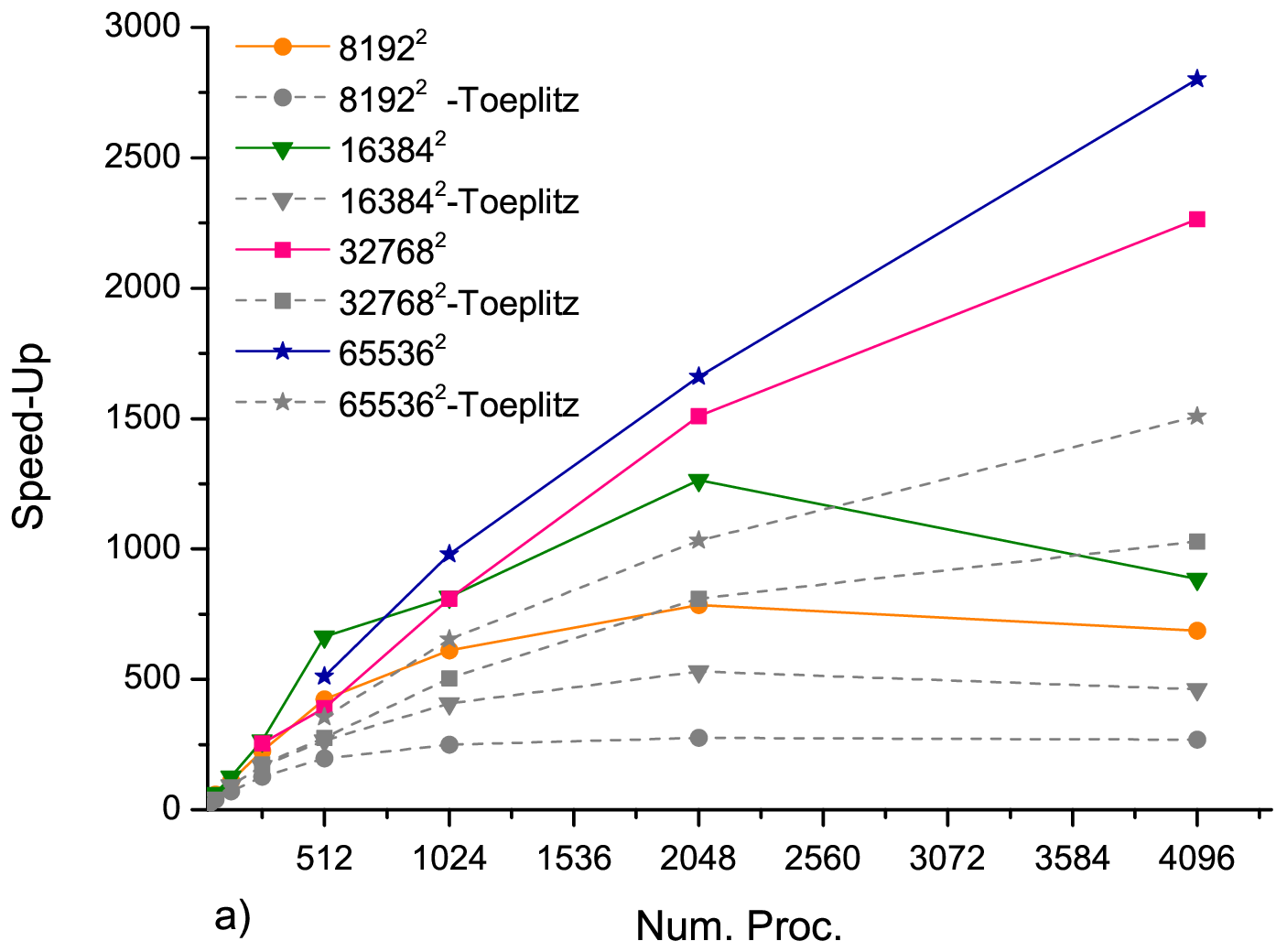} \hfill
\includegraphics[width=0.5\textwidth]{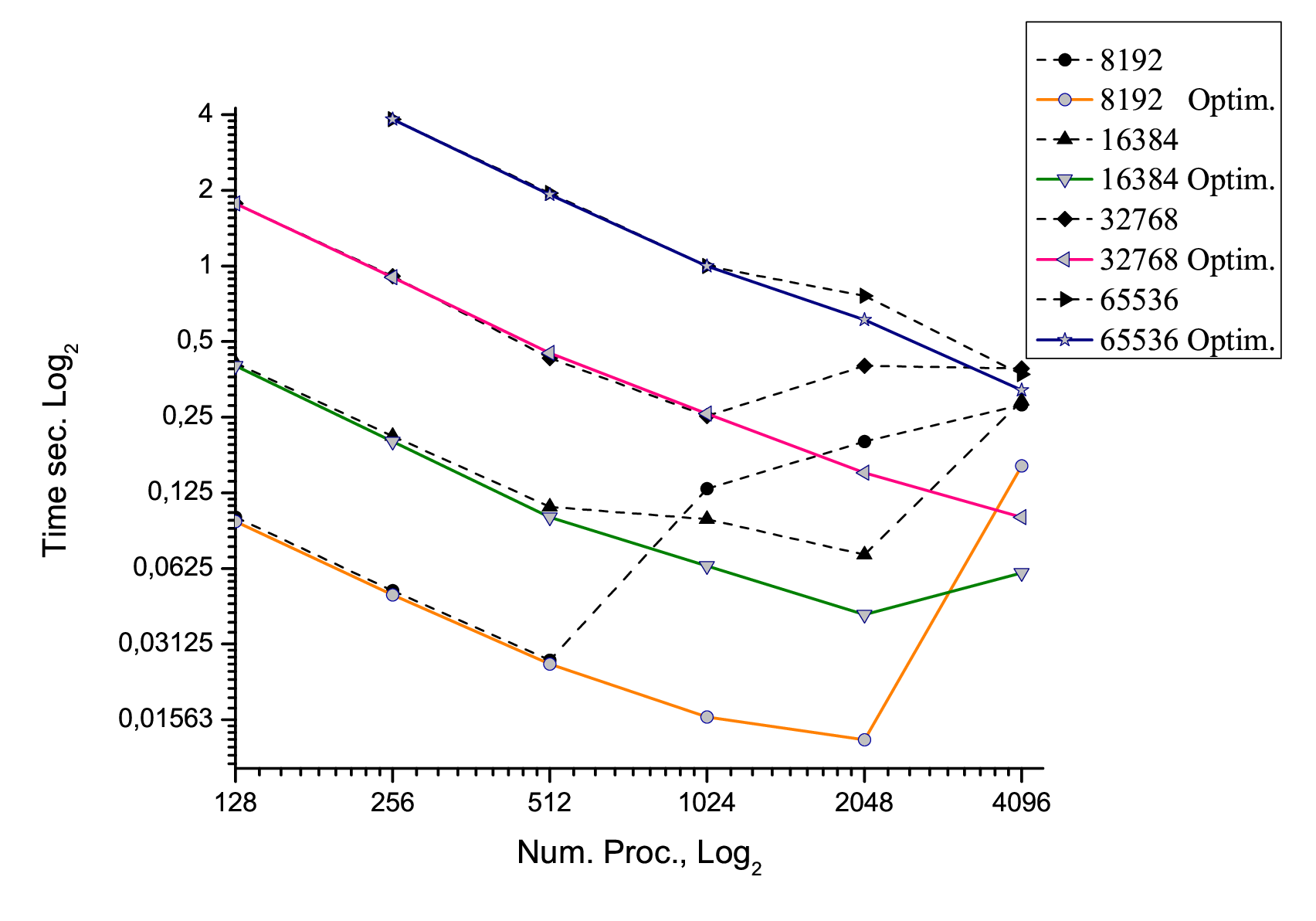}\\
\caption{a) Dependence of the speedup  on the number of the
processors for the 2D Variable-Separation Method for general (without the preliminary stage time) and
Toeplitz (including the preliminary stage time) matrices. b) The influence of the dynamic optimization of
inter-processor interactions on the computation time.} \label{main_pic2}
\end{figure}

For the 2d problem the dependencies of the time of the preliminary computations, the
time of the dichotomy process
($\mathbf{T}^{\mathrm{General}}_{\mathrm{Stage_1}},\,\mathbf{T}^{\mathrm{Toeplitz}}_{\mathrm{Stage_1}}$,
$\mathbf{T}_{\mathrm{Stage_2}}$), and the speedup
%($\mathbf{S}_{\mathrm{Step_2}},\,\mathbf{S}^{\mathrm{Toeplitz}}_{\mathrm{Step_1+Step_2}}$)
on the number of the processors are given in
Tables~\ref{table1} and in Fig.~\ref{main_pic2}. We
can see from Tables~\ref{table1} that the time ($\mathbf{T}^{\mathrm{General}}_{\mathrm{Stage_1}}$)
necessary for the computation of the vectors
$\mathbf{G}^{\mathrm{R,L}}_m,\, \mathbf{Z}^{\mathrm{R,L}}_m$  for
general matrices is independent of the number of processors and is
proportional to the number of unknowns. This is due to the fact
that the total dimension of the auxiliary vectors has order
$O(N)$. Thus, if we neglect the fact that the
matrices are Toeplitz for (\ref{poisson_test}) then the time costs
for the preliminary computations will be $O(N_1N_2)$. For Toeplitz
matrices, the use of (\ref{Zdiff}),(\ref{Greendiff3}), in comparison with the general case
(\ref{green_v}),(\ref{z_vector}) makes possible to decrease the
computation time of the auxiliary vectors
($\mathbf{T}^{\mathrm{General}}_{\mathrm{Stage_1}}$ vs.
$\mathbf{T}^{\mathrm{Toeplitz}}_{\mathrm{Stage_1}}$) by several
orders. As a result, it becomes possible to efficiently solve SLAEs with both one and several right-hand sides because the time needed for the preliminary stage of the dichotomy algorithm for the Toeplitz matrices becomes comparable to that needed for carrying out the dichotomy process(Fig.~\ref{pic:example_dichotomy2}) for one right-hand side. For solving a SLAE with one right-hand side, we should take into account the preliminary dichotomy algorithm costs. In this case, the speedup value will be  $1.5-2.5$  times less then that in solving a SLAE for several right-hand sides, i.e. when the preliminary computations independent of the number of the right-hand-sides can be neglected.

\begin{table}[!h]
\center \small
\begin{tabular}{l|cc|cc|cc|cc}
  \hline
$N_1\times N_2$ & \multicolumn{2}{c|}{8192x8192}& \multicolumn{2}{c|}{16384x16384}&\multicolumn{2}{c|}{32768x32768}&\multicolumn{2}{c}{65536x65536}\\ \hline
$\mathbf{T}^{\mathrm{General}}_{\mathrm{Stage_1}}$ &\multicolumn{2}{c|}{ $\approx3.3$sec.}&\multicolumn{2}{|c}{$\approx 13$ sec.}&\multicolumn{2}{|c}{$\approx 53$ sec.}&\multicolumn{2}{|c}{$\approx 246$ sec.}\\ \hline num. proc. &  $
\mathbf{T}^{\mathrm{Toeplitz}}_{\mathrm{Stage_1}}$&  $
\mathbf{T}_{\mathrm{Stage_2}}$  &$\mathbf{T}^{\mathrm{Toeplitz}}_{\mathrm{Stage_1}}$&$
\mathbf{T}_{\mathrm{Stage_2}}$ & $
\mathbf{T}^{\mathrm{Toeplitz}}_{\mathrm{Stage_1}}$&
 $ \mathbf{T}_{\mathrm{Stage_2}}$  & $
\mathbf{T}^{\mathrm{Toeplitz}}_{\mathrm{Stage_1}}$& $ \mathbf{T}_{\mathrm{Stage_2}}$ \\
  \hline
  32 & 0.2&3.9e-01&0.75&1.6&-&-&-&-\\
  64 & 0.1& 1.9e-01&0.33&0.84&-&-&-&-\\
  128 &5.6e-02& 1.0e-01&0.17&0.4&0.81&1.77&-&-\\
  256 &3.7e-02&4.9e-02&0.12&0.2&0.41&0.89&1.62&3.83\\
  512 & 3.0e-02&2.6e-02&8.0e-02&0.12&0.24&0.58&0.84&1.92\\
  1024 &2.6e-02&1.6e-02&6.5e-02&6.5e-02&0.17&0.28&0.51&1.0\\
  2048 &2.6e-02&1.4e-02&5.8e-02&4.2e-02 &0.13&0.15&0.36&0.59\\
  4096 &2.5e-02&1.6e-02&5.5e-02&6.0e-02&0.12&0.10&0.30&0.35\\
  \hline
\end{tabular}
\caption{The dependence of the computation time $\mathbf{T}$ on the number of processors for the 2D Poisson equation. Here
$\mathbf{T}^{\mathrm{Toeplitz}}_{\mathrm{Stage_1}}$,
$\mathbf{T}^{\mathrm{General}}_{\mathrm{Stage_1}}$ are the times of
the preliminary stage of the dichotomy algorithm for the Toeplitz matrices
and general matrices respectively; $\mathbf{T}_{\mathrm{Stage_2}}$
is the time of the implementation of the dichotomy process for any
tridiagonal matrix.}
 \label{table1}
\end{table}

We see that the
dichotomy algorithm ensures a speedup close to that linear. For
a large number of processors, no substantial decrease of
efficiency occurs because of the predominance of inter-processor
exchanges. This is ensured by dynamic optimization of the
communication interactions. Indeed, calling the function
MPI\_Reduce for the first time, we can collect information about
the time of the interaction between the processors, the volume of
the data transferred etc. for each communicator. This will make possible to optimize the processes of data exchange for subsequent calls of
MPI\_Reduce \cite{Collective:1,Tuning:1,Tuning:2,Tuning:3,CPE:CPE1206}.
The dependence of the computation time on the number of processors
in the case of the use of dynamic optimization and without it is
given in Fig.~\ref{main_pic2}.b. For a small number of processors
$p<512$, dynamic optimization does not influence the computation
time much since, in this case, the computational costs exceed
those communicational. However, as the number of the processors
increases, the effect caused by optimization becomes more than
substantial. Thus, the high efficiency of the dichotomy algorithm is provided,
on the one hand, by the low costs of the synchronization of the
computations, and on the other, by the possibility of reduction
for the data transfer time by static and dynamic optimization of
the communication interactions.

\begin{figure}[!h]
\includegraphics[width=0.5\textwidth]{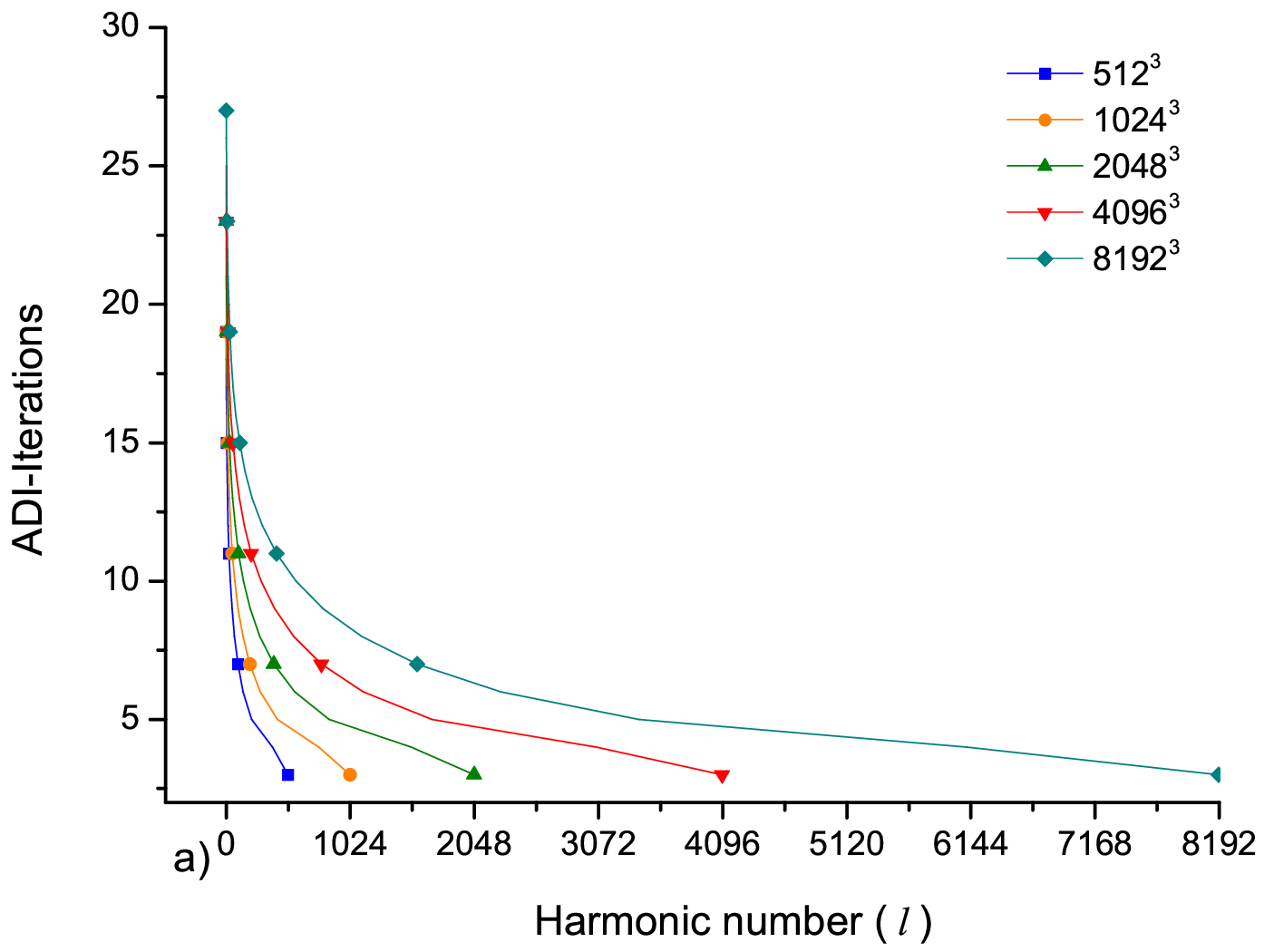}\hfill
\includegraphics[width=0.5\textwidth]{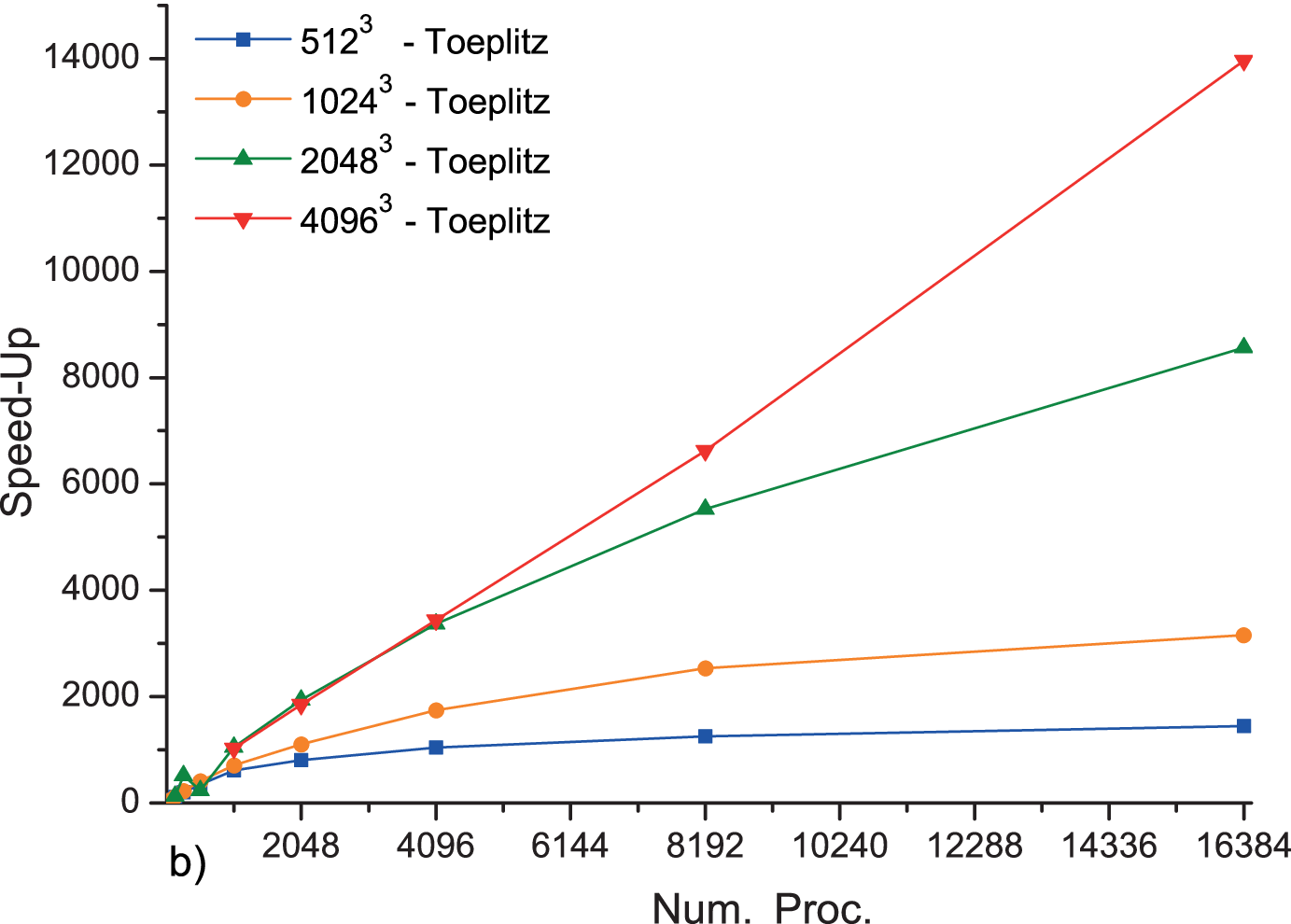} \\
%\parbox{0.5\textwidth}{\center a)}\hfill
%\parbox{0.5\textwidth}{\center b)}
\caption{a)	The number of ADI method iterations versus the harmonic number in solving problems of the form (\ref{razdiffeq22}).
b) Speedup versus the number of processors in solving the 3D Poisson equation.} \label{pic:3d}
\end{figure}

For the case of solving the 3D Poisson equation, the dependence of the calculating time and the speedup on the number of processors is given in Table~\ref{table2} and in Fig.~\ref{pic:3d}.b.  As for the 2D problem, the speedup nearly linearly depends on the number of processors. Owing to the high scalability of the dichotomy algorithm, the calculation time, with increasing number of processors, reaches a stationary value and then does not grow. The use of the efficient modification of the preliminary procedure of the dichotomy algorithm for the Toeplitz tridiagonal matrices has made it possible to reduce the time of preliminary calculations from several hours to several fractions of second.
\begin{table}[!h]
\center \small
\begin{tabular}{l|cc|cc|cc|cc|cc}
  \hline
$N_1\times N_2\times N_3$ & \multicolumn{2}{c|}{$512^3$}& \multicolumn{2}{c|}{$1024^3$}&\multicolumn{2}{c|}{$2048^3$}&\multicolumn{2}{c|}{$4096^3$}&\multicolumn{2}{c}{$8192^3$}\\ \hline
 num. proc.&  $
\mathbf{T}^{\mathrm{Toeplitz}}_{\mathrm{Stage_1}}$&  $
\mathbf{T}_{\mathrm{Stage_2}}$  &$\mathbf{T}^{\mathrm{Toeplitz}}_{\mathrm{Stage_1}}$&$
\mathbf{T}_{\mathrm{Stage_2}}$ & $
\mathbf{T}^{\mathrm{Toeplitz}}_{\mathrm{Stage_1}}$&
 $ \mathbf{T}_{\mathrm{Stage_2}}$  & $
\mathbf{T}^{\mathrm{Toeplitz}}_{\mathrm{Stage_1}}$& $ \mathbf{T}_{\mathrm{Stage_2}}$&$
\mathbf{T}^{\mathrm{Toeplitz}}_{\mathrm{Stage_1}}$&$ \mathbf{T}_{\mathrm{Stage_2}}$ \\
  \hline
  64   &2.1e-2&0.76   &6.5e-2&4.8&-&-&-&-&-&-\\
  128  &1.9e-2&0.42   &5.5e-2&2.5&1.6e-1&31.1&-&-&-&-\\
  256  &1.8e-2&0.23   &4.6e-2&1.3&1.3e-1&7.7&-&-&-&-\\
  512  &1.7e-2&0.14   &4.2e-2&0.75&1.1e-1&16.6&-&-&-&-\\
  1024 &1.5e-2&0.08   &3.8e-2&0.44&9.5e-2&3.8&2.3e-1 &24&-&-\\
  2048 &1.5e-2&0.06   &3.7e-2&0.28&8.7e-2&2.0&2.2e-1&13.3&-&-\\
  4096 &1.6e-2&4.6e-2 &4.1e-2&0.17&8.6e-2&1.1&2.1e-1&7.1&-&-\\
  8192 &1.6e-2&3.9e-2 &3.3e-2&0.12&8.5e-2&0.72&2.1e-1&3.71&-&-\\
  16384&1.6e-2&3.3e-2 &3.7e-2&9.8e-2&8.3e-2&0.46&1.8e-1&1.76&0.37&23.9\\
  \hline
\end{tabular}
\caption{The dependence of the computation time $\mathbf{T}$ on the number of processors for solving of the 3D Poisson equation. Here
$\mathbf{T}^{\mathrm{Toeplitz}}_{\mathrm{Stage_1}}$ is the time of
the preliminary stage of the dichotomy for the Toeplitz matrices; $\mathbf{T}_{\mathrm{Stage_2}}$
is the time of the implementation of the dichotomy process for any
tridiagonal matrix.}
 \label{table2}
\end{table}

\section{Conclusions}
In this article, we have proposed  a parallel algorithm for
solving tridiagonal systems of equations with Toeplitz matrices
demonstrating high efficiency including for thousands of
processors. The method is based on the Dichotomy Algorithm devised
for solving a series of tridiagonal systems of linear equations
with constant matrix and different right-hand sides. The fact
that, for Toeplitz matrices, each component of the solution vector
may be calculated without solving the SLAE makes it possible to
substantially reduce  the computation time at the preliminary step
of the Dichotomy Algorithm. As a result, it became possible to
effectively solve not only a series but a single system of
equations. The reduction of the time for the preliminary
computations is also very important in solving  a series of
problems with a large number of unknowns, since, if we disregard
the special structure of the matrix, the delay before the direct
start of solving  the equation may be substantial.

In some algorithms considered above, the presence of diagonal
dominance for the matrix of the SLAE is taken into account, which
enables us to reduce the number of interprocessor exchanges.
Sometimes this imposes the constraint that the size if systems of
the form of (\ref{systemb}) must be at least some threshold value;
otherwise, we have an accuracy loss \cite{Tridiagonal3}. Thus, the
maximal number of the processors that may be used for solving the
problem depends on the presence of diagonal dominance. If we use
the Dichotomy Algorithm, such dependencies do not arise, and the
solution satisfies (with the computer accuracy in mind) the
initial system of equations for any number of the processors. As a
result, the process of paralleling of the already existing
numerical methods including the inversion procedure of a
tridiagonal SLAE is substantially simplified.

Computational experiments confirm that dynamic optimization on the
level of the MPI-library makes it possible to substantially reduce
the time of interprocessor exchanges. The effect caused by
optimization is especially noticeable in computations with the use
of a large number of processors. Thus, the above-proposed
algorithm for solving tridiagonal systems of equations with
Toeplitz matrices makes it possible to attain a high computation
speed in a wide range of processors and guarantee high
transferability of the software.

\newpage
\bibliography{base}

\end{document}